\documentstyle[12pt,amsfonts,amscd]{amsart}
\newtheorem{Theorem}{Theorem}[section]
\newtheorem{Definition}[Theorem]{Definition}
\newtheorem{Proposition}[Theorem]{Proposition}

\newtheorem{Lemma}[Theorem]{Lemma}
\newtheorem{Corollary}[Theorem]{Corollary}
\theoremstyle{remark}

\newtheorem{Example}[Theorem]{Example}

\def\eps{\varepsilon}

\def\ovr{\overline}

\def\al{\alpha}
\def\gm{\gamma}

\def\bd{\partial}
\def\lm{\lambda}

\def\sm{\setminus}
\def\sbs{\subset}

\def\Sat{\operatorname{Sat}}

\def\re{{\mathbf {Re\,}}}

\def\be{\begin{enumerate}}
\def\ee{\end{enumerate}}
\def\bT{\begin{Theorem}}
\def\eT{\end{Theorem}}
\def\bP{\begin{Proposition}}
\def\eP{\end{Proposition}}
\def\bD{\begin{Definition}}
\def\eD{\end{Definition}}
\def\bE{\begin{Example}}
\def\eE{\end{Example}}
\def\bL{\begin{Lemma}}
\def\eL{\end{Lemma}}
\def\bC{\begin{Corollary}}
\def\eC{\end{Corollary}}

\begin{document}
\title{Functions holomorphic along holomorphic vector fields}%
\author{Kang-Tae Kim, Evgeny Poletsky and Gerd Schmalz}%
\address{Department of Mathematics, POSTECH, Pohang, 790-784, Korea}%
\email{kimkt@@postech.edu}%
\address{Department of Mathematics, Syracuse University,
Syracuse, NY, 13244, USA}
\email{eapolets@@syr.edu}
\address{University of New England, Mathematics, Armidale NSW 2351, Australia}
\email{gerd@@mcs.une.edu.au}%
\thanks{K.T.Kim and G.Schmalz were supported by the Scientific
visits to Korea program of the AAS and KOSEF. E. Poletsky was
supported by NSF Grant DMS-0500880. G. Schmalz gratefully
acknowledges support and hospitality of the Max-Planck-Institut
f\"ur Mathematik Bonn.}
\subjclass{32A99, 35A20, 35F99}%
\keywords{Holomorphic functions, holomorphic vector fields}
%\date{}
%\dedicatory{}%
%\commby{}%
% ----------------------------------------------------------------
\begin{abstract}
The main result of the paper is the following generalization of
Forelli's theorem \cite{F}: Suppose $F$ is a holomorphic vector
field with singular point at $p$,  such that $F$ is
linearizable at $p$ and the matrix is diagonalizable with the
eigenvalues whose ratios are positive reals. Then any function
$\phi$ that has an asymptotic Taylor expansion at $p$ and is
holomorphic along the complex integral curves of $F$ is
holomorphic in a neighborhood of $p$.
\par We also present an example to show that the requirement
for  ratios of the eigenvalues to be positive reals is
necessary.
\end{abstract}
\maketitle
\section{Introduction}
\par Let $F$ be a holomorphic vector field on a complex manifold
$M$. Consider a function $\phi$ on $M$ that is holomorphic on
the complex integral curves of $F$. In general, there is not
much to say about such functions. However, Forelli noted in
\cite{F} that if $M={\Bbb D}^N$ is the unit polydisk in ${\Bbb
C}^N$ and $F(z)=z$ is the holomorphic Euler vector field, then
$\phi$ is holomorphic on ${\Bbb D}^N$ provided $\phi$ is
infinitely smooth in a neighbourhood of the origin. In fact,
(see \cite{KW}) $\phi$ only needs to have a Taylor series at
the origin and this is essential. Counterexamples are easy to
construct. Recently, Chirka in \cite{Ch} showed that in ${\Bbb
C}^2$ the straight lines, which are the integral curves of
$F(z)=z$, can be replaced by families of transversal
holomorphic curves passing through the origin.
\par In this paper we study the question what properties
of a holomorphic vector field $F$ with a critical point $z_0$
imply a Forelli type result. We prove that if $F$ is
linearizable at $z_0$ and if the linearization matrix is
diagonalizable and the pairwise ratios of its eigenvalues are
positive, then any function that is holomorphic on the complex
integral curves of $F$ and has a Taylor series at $z_0$ is
holomorphic in the open subset of $M$ saturated by $z_0$ (see
Section \ref{S:fhf} for definitions). Note that the integral
curves of generic vector field do not pass through the origin
and Chirka's result does not apply.
\par An example, presented in Section \ref{S:ce}, of a holomorphic
vector field and a non-holomorphic $C^\infty$-function
holomorphic along this field shows that the imposed
restrictions on the eigenvalues are necessary.
\par The proof follows the original steps of Forelli. The main
difficulty to overcome is the absence of a multitude of
separatrices in the general case. So we had to introduce the
appropriate notion of an asymptotic expansion. This notion
replaces the Taylor series of holomorphic functions on
separatrices. An elementary theory of such expansions is
contained in Sections \ref{S:ae} and \ref{S:aehf}. The main
outcome of this theory is the estimate of the coefficients in
Theorem \ref{T:ec} absolutely similar to the classical Cauchy
inequalities for holomorphic functions on the disk.
\section{Asymptotic expansions}\label{S:ae}
\bD Let $f$ be a function on the right half-plane ${\Bbb H}$.
Suppose that $\{\lm_j\}$, $j\ge1$, is a strictly increasing
sequence of positive numbers converging to infinity, $\lm_0=0$
and let $\{n_j\}$, $j\ge1$, be a sequence of positive integers.
A formal series
\begin{equation}\label{e:gas}
\sum_{j=0}^\infty \sum_{k=0}^{n_j}p_{jk}e^{-\mu_{jk}z-\nu_{jk}\ovr
z}\end{equation} is called an asymptotic expansion of $f$ if all
$\mu_{jk},\nu_{jk}\ge0$,
$\mu_{jk}+\nu_{jk}=\lm_j$ for every $j$ and $k$ and for every $n$
we have
$$\left|f(z)-\sum_{j=0}^n \sum_{k=0}^{n_j}p_{jk}e^{-\mu_{jk}z-\nu_{jk}\ovr
z}\right|e^{\lm_n\re z}\to0$$ as $\re z\to\infty$. \eD
\par The same definition can be trivially rephrased as follows:
\bP A function $f$ on ${\Bbb H}$ has an asymptotic expansion
(\ref{e:gas}) if and only if for every $\eps>0$
$$\left|f(z)-\sum_{j=0}^n \sum_{k=0}^{n_j}p_{jk}e^{-\mu_{jk}z-\nu_{jk}\ovr
z}\right|e^{(\lm_{n+1}-\eps)\re z}\to0$$ as $\re z\to\infty$ in
${\Bbb H}$.
\eP
\begin{pf} If
$$\left|f(z)-\sum_{j=0}^{n+1}\sum_{k=0}^{n_j}p_{jk}e^{-\mu_{jk}z-\nu_{jk}\ovr
z}\right|e^{\lm_{n+1}\re z}\to0,$$ then
$$\left|f(z)-\sum_{j=0}^n \sum_{k=0}^{n_j}p_{jk}e^{-\mu_{jk}z-\nu_{jk}\ovr
z}\right|e^{(\lm_{n+1}-\eps)\re z}\to0$$ as $\re z\to\infty$.
\par If
$$\left|f(z)-\sum_{j=0}^n \sum_{k=0}^{n_j}p_{jk}e^{-\mu_{jk}z-\nu_{jk}\ovr
z}\right|e^{(\lm_{n+1}-\eps)\re z}\to0$$ as $\re z\to\infty$, then
$$\left|f(z)-\sum_{j=0}^n \sum_{k=0}^{n_j}p_{jk}e^{-\mu_{jk}z-\nu_{jk}\ovr
z}\right|e^{\lm_{n}\re z}\to0.$$
\end{pf}
\par We will say that two asymptotic expansions are equal if all their
non-zero terms coincide. Let us state some simple properties of
asymptotic expansions.
\bP\label{P:spas} 1) Every function $f$ on ${\Bbb H}$ has at most
one asymptotic expansion (\ref{e:gas}).
\par 2) If a holomorphic function on ${\Bbb H}$ has an asymptotic
expansion (\ref{e:gas}), then all $\nu_{jk}=0$.
\eP
\begin{pf} 1) Suppose that $f$ has two asymptotic expansions:
$$\sum_{j=0}^\infty \sum_{k=0}^{n_j}p_{jk}e^{-\mu_{jk}z-\nu_{jk}\ovr
z}\text{ and }\sum_{j=0}^\infty
\sum_{k=0}^{m_j}q_{jk}e^{-\al_{jk}z-\beta_{jk}\ovr z},$$ where
$\mu_{jk}+\nu_{jk}=\lm_j$ and $\al_{jk}+\beta_{jk}=\gm_j$. Let
$l$ be the first number when two non-zero terms
$$\sum_{k=0}^{n_l}p_{lk}e^{-\mu_{lk}z-\nu_{lk}\ovr
z}\text{ and } \sum_{k=0}^{m_l}q_{lk}e^{-\al_{lk}z-\beta_{lk}\ovr
z}$$ do not coincide.
\par Suppose that $\gm_l>\lm_l$ and $n$ is the least integer such that
$\lm_n\ge\gm_l$. Since
$$\left|f-\sum_{j=0}^l\sum_{k=0}^{m_j}q_{jk}e^{-\al_{jk}z-\beta_{jk}\ovr z}
\right|e^{\gm_l\re z}\to 0$$ and
$$\left|f-\sum_{j=0}^n \sum_{k=0}^{n_j}p_{jk}e^{-\mu_{jk}z-\nu_{jk}\ovr
z}\right|e^{\gm_l\re z}\to 0,$$ we see that
$$\left|\sum_{j=l}^n\sum_{k=0}^{n_j}p_{jk}e^{-\mu_{jk}z-\nu_{jk}\ovr
z}-\sum_{k=0}^{m_l}q_{lk}e^{-\al_{lk}z-\beta_{lk}\ovr
z}\right|e^{\gm_l\re z}\to 0$$ as $\re z\to\infty$. Hence
$p_{jk}=0$, $l\le j\le n-1$. In particular, all $p_{lk}=0$ and
this contradicts our assumption that some of $p_{lk}\ne0$.
\par Thus $\lm_l=\gm_l$ and, letting $z=x+iy$, we get
$$\left|\sum_{k=0}^{n_l}p_{lk}e^{-i(\mu_{lk}-\nu_{lk})y}-
\sum_{k=0}^{m_l}q_{lk}e^{-i(\al_{lk}-\beta_{lk})y}\right|\to 0.$$
But this means that $n_l=m_l$ and $p_{lk}=q_{lk}$.
\par 2) Let $m$ be the least integer such that $\nu_{mk}\ne0$ for some $k$.
After subtracting from $f$ the first $n$ holomorphic terms of
its asymptotic expansion we may assume that
$$f(z)=\sum_{k=0}^{n_m}p_{mk}e^{-\mu_{mk}z-\nu_{mk}\ovr z}+g(z),$$ where
all $\nu_{mk}>0$ and $g(z)e^{(\lm_m+\eps)\re z}\to0$ as $\re
z\to\infty$ for some $\eps>0$.  The integral of the function
$$h(z)=f(z)e^{\lm_mz}=
\sum_{k=0}^{n_m}p_{mk}e^{2i\nu_{mk}y}+g(z)e^{\lm_mz}$$ by $dz$
over the boundary of the rectangle $\{t-1\le x\le t+1, -s\le
y\le s\}$, where $t>1$ and $s>0$, is equal to
$$-4i\sum_{k=0}^{n_m}p_{mk}\sin2\nu_{mk}s+o(e^{-\eps t})=0.$$
But this is possible if and only if all $\nu_{mk}=0$.
\end{pf}
\section{Asymptotic expansions of holomorphic functions}\label{S:aehf}
\par As it follows from the previous section an asymptotic expansion
of a holomorphic function $f$ on ${\Bbb H}$ is a formal series
\begin{equation}\label{e:chas}
\sum_{j=0}^\infty c_je^{-\lm_jz}.
\end{equation}
\par We will need a simple version of the maximum principle.
\bL\label{L:vmp} Suppose that $f$ is a bounded holomorphic function
on ${\Bbb H}$, $\lm>0$, and
$$\limsup_{\re z\to\infty}f(z)e^{\lm\re z}<\infty.$$
Let $f^*(iy)$ be the non-tangential boundary values of $f$ on
$i{\Bbb R}$. Then $|f(z)|\le Me^{-\lm\re z}$ on ${\Bbb H}$,
where $M=\sup_{y\in{\Bbb R}}|f^*(iy)|$.
\eL
\begin{pf} Note that the function $g(z)=f(z)e^{\lm z}$ is bounded on
${\Bbb H}$ and $|g^*(iy)|\le M$, $y\in{\Bbb R}$. Hence
$|g(z)|\le M$ on ${\Bbb H}$.
\end{pf}
\bC\label{C:aei} Let $f$ be a bounded holomorphic function on
${\Bbb H}$ with an asymptotic expansion $\sum_{j=0}^\infty
c_je^{-\lm_jz}$. Let $f_n(z)=f(z)-\sum_{j=0}^nc_je^{-\lm_jz}$.
Then $|f_n(z)|\le M_ne^{-\lm_{n+1}\re z}$ on ${\Bbb H}$, where
$M_n=\sup_{y\in{\Bbb R}}|f_n^*(iy)|$.
\eC
\begin{pf} By Lemma \ref{L:vmp} for every $\eps>0$ we have
$|f_n(z)|\le M_ne^{-(\lm_{n+1}-\eps)\re z}$ on ${\Bbb H}$.
Letting $\eps$ go to 0 we get our result
\end{pf}
\bC\label{C:uae} If an asymptotic expansion of $f$ has zero
coefficients, then $f\equiv0$ on ${\Bbb H}$.
\eC
\par The theorem below gives us estimates for the coefficients. We start
with a simple lemma.
\bL\label{L:it} Let $f$, $|f|\le M$, be a holomorphic function on
${\Bbb H}$ with an an asymptotic expansion $\sum_{j=0}^\infty
c_je^{-\lm_jz}$. If $a$ is real then the function
$$g(z)=f(z+ia)-f(z)$$ is holomorphic on ${\Bbb H}$,
$|g|\le2M$ there, and $g$ has the asymptotic expansion
$$\sum_{j=0}^\infty(e^{-ia\lm_j}-1)c_je^{-\lm_jz}.$$
\eL
\begin{pf} Clearly, the function $g$ is holomorphic and $|g|\le2M$ on
${\Bbb H}$. Note that
$$\sum_{j=0}^k
c_je^{-\lm_j(z+ia)}-\sum_{j=0}^k
c_je^{-\lm_jz}=\sum_{j=0}^k(e^{-ia\lm_j}-1)c_je^{-\lm_jz}.$$
Hence by Corollary \ref{C:aei} we get
\begin{equation}\begin{align}
&\left|g(z)-\sum_{j=0}^k(e^{-ia\lm_j}-1)c_je^{-\lm_jz}\right|\notag\\&
\le\left|f(z)-\sum_{j=0}^k
c_je^{-\lm_jz}\right|+\left|f(z+ia)-\sum_{j=0}^k
c_je^{-\lm_j(z+ia)}\right|\le Ce^{-\lm_{k+1}\re z}.\notag
\end{align}\end{equation}
\end{pf}
\bT\label{T:ec} Let $f$, $|f|<M$, be a holomorphic function on ${\Bbb H}$
with an asymptotic expansion $\sum_{j=0}^\infty
c_je^{-\lm_jz}$. Then $|c_j|\le M$.
\eT
\begin{pf} The proof is by induction on the indices of the
coefficients. Clearly $|c_0|\le M$. Suppose the theorem is
proved for all $c_j$, $j\le m$, for all functions $f$ and
constants $M$ as set out in the theorem. Let $a=\pi/\lm_{m+1}$.
Then $e^{-ia\lm_{m+1}}-1=-2$. If
$$g(z)=f(z+ia)-f(z)$$ and $\sum_{j=0}^\infty g_je^{-\lm_jz}$ is
the asymptotic expansion of $g$, then by Lemma \ref{L:it}
$g_0=0$. Hence, by Lemma \ref{L:vmp}, the function
$h(z)=e^{\lm_1z}g(z)$ is bounded on ${\Bbb H}$ and $|h|\le2M$
on ${\Bbb H}$. The asymptotic expansion $\sum_{j=0}^\infty
h_je^{-\mu_jz}$ of $h$ is equal to $\sum_{j=1}^\infty
g_je^{-(\lm_j-\lm_1)z}$. Thus
$$h_m=g_{m+1}=-2c_{m+1}.$$ By the induction assumption
$$|c_{m+1}|\le\frac{|h_m|}2\le M.$$
\end{pf}
\bC\label{C:ct} Suppose that a bounded holomorphic
function $f$ on ${\Bbb H}$ has an asymptotic expansion
$\sum_{j=0}^\infty c_je^{-\lm_jz}$ such that the series
$\sum_{j=0}^\infty e^{-\lm_jd}$ converges for some $d>0$. Then
the series $\sum_{j=0}^\infty c_je^{-\lm_jz}$ converges
uniformly to $f(z)$ on ${\Bbb H}_d=\{\re z\ge d\}$.
\eC
\section{Taylor series and asymptotic expansions}
\par Let $\phi$ be a function defined in a neighborhood of the
origin in ${\Bbb C}^N$. We say that $\phi$ has a Taylor series
at the origin if there is a formal series
\begin{equation}\label{e:ts}
\sum_{j=0}^\infty\sum_{|k|+|m|=j}a_{km}z^k\ovr z^m\end{equation}
such that for every $j$
$$\left|\phi(z)-\sum_{j=0}^n\sum_{|k|+|m|=j}a_{km}z^k\ovr z^m\right|=
o(|z|^n).$$ (Here for $k=(k_1,\dots,k_N)$ and
$z=(z_1,\dots,z_N)$, $z^k=z_1^{k_1}\cdots z_N^{k_N}$.)
\par Let $\alpha=(\alpha_1,\dots,\alpha_N)$ be a vector in $\mathbb R^N$
with positive components and $c=(c_1,\dots,c_N)$ a vector in
$\mathbb D^N$. We define the function $s_c:\quad \mathbb C\to
\mathbb C^N$ by
$$s_c: \quad  \zeta \mapsto (c_1e^{-\alpha_1
\zeta},\dots,c_Ne^{-\alpha_N \zeta}).$$ Notice that
$s_c(\zeta)\in {\mathbb D}^N$ if $\zeta\in\mathbb H$.
\par Denote by $\{\lambda_j\}$ the sequence of all possible
values
$$(\alpha,k)+ (\alpha,m)= \alpha_1k_1+\cdots+ \alpha_Nk_N+
\alpha_1m_1+\cdots+ \alpha_Nm_N,$$
where $k,m$ are multi-indices, in ascending order.
\bP\label{P:aets} Let $\phi$ be a function defined on $\mathbb
D^N$. If $\phi$ has the Taylor series (\ref{e:ts}) at the
origin, then the function $f_c= \phi\circ s_c:\quad \mathbb
H\to \mathbb C$ has the asymptotic expansion
$$\sum_{j=0}^\infty\left(\sum_{(\al,k)+(\al,m)=\lm_j}a_{km}c^k\ovr
c^m e^{-(\al,k)\zeta-(\al,m)\ovr\zeta}\right)$$
on ${\Bbb H}$.
\eP
\begin{pf} Without any loss of generality we may assume that $\al_1$ is
the minimal number among those $\al_j$ for which $c_j\ne0$.
Then $|s_c(\zeta)|\sim e^{-\al_1\re\zeta}$.
\par For any $p$ choose $n$ such that $n\al_1>\lm_p$. By the definition
of the Taylor series
$$\left(f_c(\zeta)-\sum_{j=0}^n\left(\sum_{|k|+|m|=j}a_{km}c^k\ovr
c^m
e^{-(\al,k)\zeta-(\al,m)\ovr\zeta}\right)\right)e^{n\al_1\re\zeta}\to0$$
as $\re\zeta\to\infty$.
\par Note that if $|k|+|m|>n$ then $(\al,k)+(\al,m)>na_1>\lm_p$. Hence all
multi-indices $k$ and $m$ for which $(\al,k)+(\al,m)\le\lm_p$
are included into the sum
$$\sum_{j=0}^n\left(\sum_{|k|+|m|=j}a_{km}c^k\ovr
c^m e^{-(\al,k)\zeta-(\al,m)\ovr\zeta}\right).$$ Splitting this
sum into the sum over all $k$ and $m$ such that
$(\al,k)+(\al,m)\le\lm_p$ and the sum over all $k$ such that
$(\al,k)+(\al,m)>\lm_p$ and noticing that the terms of the
second sum multiplied by $e^{\lm_p\re\zeta}$ go to 0 as
$\re\zeta\to\infty$, we see that
$$\left(f_c(\zeta)-\sum_{(\al,k)+(\al,m)\le\lm_p}a_{km}c^k\ovr
c^m
e^{-(\al,k)\zeta-(\al,m)\ovr\zeta}\right)e^{\lm_p\re\zeta}\to0$$
as $\re\zeta\to\infty$.
\par Rearranging terms in the latter sum so that
\begin{equation}\begin{align}
&\sum_{(\al,k)+(\al,m)\le\lm_p}a_{km}c^k\ovr c^m
e^{-(\al,k)\zeta-(\al,m)\ovr\zeta}\notag\\&
=\sum_{j=0}^p\left(\sum_{(\al,k)+(\al,m)=\lm_j}a_{km}c^k\ovr c^m
e^{-(\al,k)\zeta-(\al,m)\ovr\zeta}\right)\notag
\end{align}\end{equation} we obtain the desired result.
\end{pf}
\par From the second statement in Proposition \ref{P:spas} we
immediately get
\bC\label{C:hts} If in Proposition \ref{P:aets} the function
$f_c$ is holomorphic, then
$$\sum_{(\al,m)=\lm}a_{km}{\ovr c}^m\equiv0$$ for all $\lm>0$.
\eC
\section{Linear systems of differential equation}\label{S:lsde}
\par The goal of this section is to prove the following
theorem.
\bT\label{T:hts} Let $A$ be a diagonal $N\times N$-matrix with eigenvalues
$\al_1,\dots,\al_N$ such $\al_j/\al_k>0$ for every pair $(j,k)$
and let $\phi$ be a function  defined on a neighborhood of the
origin. If $\phi$ has the Taylor series (\ref{e:ts}) at the
origin and is holomorphic on all solutions of the system
$z'=Az$, then this Taylor series does not contain
anti-holomorphic variables and converges on a neighborhood of
the origin to the function $\phi$.
\eT
\begin{pf} Changing the time variable $\zeta$ to $-\al_1\zeta$ we may
assume that all eigenvalues of $A$ are real and negative. Also
without any loss of generality we may assume that the
coordinates on ${\Bbb C}^N$ are chosen so that $\phi$ is
defined on the unit polydisk ${\Bbb D}^N$.
\par It is well-known that for every point $c\in{\Bbb C^N}$
there is a solution $s_c(\zeta):\,{\Bbb C}\to{\Bbb C}^N$ of our
differential equation which can be written in the form
$s_c(\zeta)=(c_1e^{\al_1\zeta},\dots,c_Ne^{\al_N\zeta})$.
\par By Corollary \ref{C:hts}
$$\sum_{(\al,m)=\lm}a_{km}{\ovr c}^m\equiv0$$ for all $\lm>0$
and all $c\in{\Bbb C}^N$. But $c^k$ are linearly independent
monomials over ${\Bbb C}^N$ and, therefore, the identity above
holds only if all $a_{km}=0$ when $m\ne0$.
\par To prove the second part we notice that $\phi$ is
bounded in some neighborhood of the origin since it has a
Taylor expansion.  Therefore we may assume that $|\phi|\le M$
on ${\Bbb D}^N$. By Proposition \ref{P:aets} and Theorem
\ref{T:ec} the absolute values of all polynomials
$$\sum_{(\al,k)=\lm_j}a_{k0} c^k$$
do not exceed $M$. By the Cauchy estimates $|a_{k0}|\le1$. Thus
the series $\sum_ka_{k0}z^k$ absolutely converges on compacta
in ${\Bbb D}^N$ to a function $\psi$. Hence, the formal series
$$\sum_{j=0}^\infty\left(\sum_{(\al,k)=\lm_j}a_{k0}c^k\right)
e^{-\lm_j\zeta}$$ of the asymptotic expansion for
$\phi_c(\zeta)=\phi(s_c(\zeta))$ converges on ${\Bbb H}_d$ for
every $d>0$ to a function $\psi_c(\zeta)$. But this means that the
difference $\phi_c(\zeta)-\psi_c(\zeta)$ has an asymptotic
expansion with zero coefficients. By Corollary \ref{C:uae}
$\phi_c\equiv\psi_c$ and $\phi\equiv\psi$. Thus $\phi$ is
holomorphic on ${\Bbb D}^N$.
\end{pf}
{\bf Remark.} The first part of this theorem can be proved when
there is $\al\in{\Bbb C}$ such that all numbers $\al_j\al$ have
negative real parts. As the following example shows this
requirement is rather essential. Consider the differential
equation $z_1'=z_1,\,z_2'=-z_2$ in ${\Bbb C}^2$. The solutions
are: $z_1(\zeta)=c_1e^\zeta$ and $z_2(\zeta)=c_2e^{-\zeta}$. The
function $\phi(z_1,z_2)=\ovr z_1\ovr z_2$ is real analytic,
constant on solutions, and its Taylor series is anti-holomorphic.
\section{$F$-holomorphic functions}\label{S:fhf}
\par Let $F$ be a holomorphic vector field on a complex manifold $M$ of
dimension $N$.  We say that a function $\phi$ on $M$ is
$F$-holomorphic if the restrictions of $\phi$ to the integral
curves of the $F$ are holomorphic. If $\phi$ is differentiable
then it means that its derivative along $\ovr F$ is equal to
$0$ or in local coordinates $(z_1,\dots,z_N)$ on $M$  if
$F=(F_1,\dots,F_N)$, then
$$\nabla_{\ovr F}\phi=\sum_{j=1}^N\frac{\bd\phi}{\bd\ovr
z_j}\ovr F_j\equiv0.$$
\par The maximal integral curves of $F$ form a holomorphic
foliation of $M$ (see \cite[Ch. 1]{IY}). If $U$ is an open set in
$M$ then following \cite[Def. 2.16]{IY} we denote by $\Sat(U,F)$
the union of all maximal integral curves of $F$ intersecting $U$.
By \cite[Lemma 2.17]{IY} $\Sat(U,F)$ is open.
\bL\label{L:es} Let $M$ be a complex manifold and let $F$ be
a holomorphic vector field on $M$. If an $F$-holomorphic function
$\phi$ on $M$ is holomorphic on an open set $U\sbs M$, then it is
holomorphic on $\Sat(U,F)$.
\eL
\begin{pf} Let $V\sbs \Sat(U,F)$ be the open set, where $\phi$ is
holomorphic and let $S$ be a maximal integral curve of $F$ in
$M$ intersecting $U$. If $S\cap V\ne S$ then there is a point
$z_0\in S$ lying in the relative boundary of $S\cap V$ in $S$.
If $z_0$ is a critical point of $F$, then $S=\{z_0\}\sbs U$ and
$\phi$ is holomorphic near $z_0$. If $z_0$ is not a critical
point of $F$, then by the Rectification Theorem (see
\cite[Theorem 1.18]{IY}) there are local coordinates
$(z_1,\dots,z_N)$ around $z_0$ in $M$, where
$F(z)=(1,0,\dots,0)$. We may assume that these local
coordinates are defined on ${\Bbb D}^N$ and $S\cap{\Bbb
D}^N=\{(\zeta,0,\dots,0)\}$. Since $z_0$ is a relative boundary
point of $S\cap V$ in $S$, there is $\zeta_0\in{\Bbb D}$ and an
open neighborhood $W$ of $(\zeta_0,0,\dots,0)\}$ in ${\Bbb
D}^N$, where $\phi$ is holomorphic. By Hartogs Lemma (see
\cite[Lemma 3.6.2.1]{Sh}) there is $\eps>0$ such that $\phi$ is
holomorphic on an open neighborhood of $S\cap{\Bbb D}^N$. Hence
$z_0\in V$ and, by contradiction, $S\sbs V$. Thus
$\Sat(U,F)=V$.
\end{pf}
\par In general, $F$-holomorphic functions need not to be holomorphic.
However, an additional information on their behavior near special
critical points of $F$ implies that they are holomorphic on rather
large sets.
\par Let us recall that a vector field $F$ is linearizable at a
critical point $z_0\in M$ of $F$ if there is a biholomorphic
mapping $H$ of a neighborhood $U$ of $z_0$ onto a neighborhood
$V$ of $H(z_0)$ such that in new coordinates  $F(w)=Aw$, where
$A$ is an $N\times N$-matrix.
\par If $z_0$ is a critical point of $F$ then we define the set
$\Sat(z_0, F)$ as the union of all maximal integral curves $L$
of $F$ satisfying: $L \cap U \not= \emptyset$ for any open
neighborhood $U$ of $z_0$. If $F$ is linearizable at its
isolated critical point $z_0$ and there is $\al\in{\Bbb C}$
such that  all numbers $\al_j\al$ have negative real parts,
where $\al_j$ are the eigenvalues of the linearization matrix
$A$, then there is a neighborhood $U$ of $z_0$ such that the
set $\Sat(z_0,F)=\Sat(U,F)$ and, consequently, is open.
\bT\label{T:mt} Let $F$ be a holomorphic vector field on a
complex manifold $M$. Suppose that $F$ is linearizable at its
isolated critical point $z_0$, the linearization matrix $A$
diagonal and the ratios $\al_j/\al_k$ of all eigenvalues of $A$
are positive. If an $F$-holomorphic function $\phi$ on $M$ has
a Taylor series at $z_0$, then $\phi$ is holomorphic on
$\Sat(z_0,F)$.
\eT
\begin{pf} By Theorem \ref{T:hts} $\phi$ is holomorphic on a
neighborhood $U$ of $z_0$ and by Lemma \ref{L:es} on
$\Sat(z_0,F)\sbs\Sat(U,F)$.
\end{pf}
\bC Suppose that a holomorphic vector field $F$ on a complex
manifold $M$ has only critical points $z_n$ described in the
statement of Theorem \ref{T:mt} and $M$ is the union of the
sets $\Sat(z_n,F)$. Then any $F$-holomorphic
$C^\infty$-function on $M$ is holomorphic. In particular, if
$M$ is compact then any such function is constant.
\eC
\section{Two counterexamples}\label{S:ce}
\par Consider the holomorphic vector field $F(z_1,z_2)=(\al z_1,\beta z_2)$
in ${\Bbb C}^2$, where $\al$ and $\beta$ are both different
from 0 and $\al/\beta$ is not a positive number. First we
consider the case $\beta=-t\alpha$ with $t>0$.  The integral
curves of $F$ do not change if we replace the complex time
$\zeta$ by $\al\zeta$. So we may assume that
$F(z_1,z_2)=(z_1,-tz_2)$. The function $|z_1|^t|z_2|$ is
constant on integral curves of $F$ and smooth when
$z_1z_2\ne0$. To make it smooth everywhere we use
$$\phi(z_1,z_2)=\exp\left(-\frac1{|z_1|^t|z_2|}\right).$$ Since
$\phi$ is constant on the integral curves it is holomorphic
along $F$.
\par If $\al/\beta$ is not real then replacing the complex time
$\zeta$ by $\tau\zeta$ with an appropriate $\tau$ we may assume
that $\al=\al_1+i\al_2$, $\al_1,\al_2\in{\Bbb R}$, $\al_1<0$,
$\al_2>0$ and $\beta=t\ovr\al$, $t>0$.
\par Let $$\gm=\frac1{2\al_1}-\frac i{2\al_2}.$$ Direct calculations
show that
\begin{equation}\label{e:tl}
\gamma\re\alpha\zeta+ \frac{\ovr\gm}{t}\re\beta\zeta \equiv \zeta.
\end{equation}
\par The sector $S=\{\zeta\in{\Bbb C}:\,\zeta=r\gm+s\ovr\gm,
r,s<0\}$ lies in the left half-plane and has an angle strictly
less than $\pi$. Hence there is $b>1$ such that function
$\zeta^b$ is defined with $\re\zeta^b<0$ on $S$. Hence the
function
$$\phi(z_1,z_2)=\exp\left(\gm\log|z_1|+\frac{\ovr\gm}t\log|z_2|
\right)^b$$ is defined when $(z_1,z_2)\in{\Bbb D}^2$ and
$z_1z_2\ne0$. If $z_1z_2=0$ we let $\phi=0$. The function
$\phi$ is not holomorphic, but we shall demonstrate that it is
smooth and holomorphic along every integral curve F.
\par Clearly, $\phi\in C^\infty({\Bbb D}^2\sm\{z_1z_2=0\})$. But
it also is $C^\infty$ on $\{z_1z_2=0\}$. To see this we note
that if $\xi=\gm\log|z_1|+\frac{\ovr\gm}t\log|z_2|\in S$, then
$$\re \xi^b\le
-c\max\{(-\log|z_1|)^b,(-\log|z_2|)^b\},$$ where $c>0$.
Thus if, say, $z_1\to0$, then
$$|\phi(z_1,z_2)|\le e^{-c(-\log|z_1|)^b}=
|z_1|^{c(-\log|z_1|)^{b-1}},$$ i.e., $\phi$ decreases faster
than any power of $z_1$. This shows that $\phi\in
C^\infty({\Bbb D}^2)$ and its Taylor series on $\{z_1z_2=0\}$
is zero.
\par  The integral curves of $F$ admit a parametrization:
$z_1(\zeta)=C_1e^{\al\zeta}$ and
$z_2(\zeta)=C_2e^{\beta\zeta}$. By (\ref{e:tl})
$$\gm\log|z_1(\zeta)|+\frac{\ovr\gm}{t}\log|z_2(\zeta)|=
\zeta+\gm\log|C_1|+\frac{\ovr\gm}{t}\log|C_2|.$$
Hence
$$\phi(z_1(\zeta),z_2(\zeta))=
\exp\left(\zeta+\gm\log|C_1|+\frac{\ovr\gm}{t}\log|C_2|\right)^b$$
is a holomorphic function of $\zeta$.

\end{document}